\documentclass[12pt,letterpaper]{amsart}
\usepackage{epic,eepic,latexsym, amssymb, amscd, amsfonts, xypic, graphicx}
\usepackage[usenames]{color}

 \newlength{\baseunit}               
 \newcount{\numlines}                
\newcommand{\tpoint}[1]{\vspace{3mm}\par \noindent \refstepcounter{subsection}{\bf \thesubsection.}
  {\em #1. ---} }

\newcommand{\bpoint}[1]{\vspace{3mm}\par \noindent \refstepcounter{subsection}{\bf \thesubsection.}
  {\bf #1.} }

\newtheorem{tm}{Theorem}[subsection]

\newtheorem{lm}[tm]{Lemma}

\newcommand{\Spec}{\operatorname{Spec}}

\newcommand{\A}{\mathcal{A}}
\newcommand{\E}{\mathcal{E}}
\newcommand{\F}{\mathcal{F}}

\newcommand{\secretnote}[1]{}

\newcommand{\lremind}[1]{{}}

\newcommand{\Hom}{\operatorname{Hom}}

\newcommand{\Oh}{\mathcal{O}}

\newcommand{\comment}[1]{}

\newcommand{\RelSpec}[2]{\ensuremath{\textbf{\textrm{Spec}}_{#1}#2}}

\begin{document}
\pagestyle{plain}

\title{\large{$S_3$-covers of Schemes}}

\date{April 27, 2008.}

\author{Robert W. Easton}

\address{Department of Mathematics, University of Utah, Salt Lake City, Utah 84112}
\email{easton@math.utah.edu}

\begin{abstract}
We analyze flat $S_3$-covers, attempting to create structures parallel to those found in the abelian theory.  We use an initial local analysis as a guide in finding a global description.
\end{abstract}
\maketitle


\vspace{-0.2in}
{\parskip=12pt 

\section{Introduction}

Given a finite group $G$, a {\em $G$-cover} of a scheme $X$ is a scheme $Y$ together with a faithful $G$-action on $Y$ and a finite $G$-equivariant morphism $\pi:Y\to X$ which identifies $X$ with the geometric quotient $Y/G$. If one considers only schemes over a fixed algebraically closed field $k$ of characteristic prime to the order of $G$, then to each $G$-cover $\pi:Y\to X$ there is a decomposition $\pi_* \Oh_Y=\oplus_{\rho\in G^{\vee}}\F_{\rho}$, where $\F_{\rho}$ is an $\Oh_X[G]$-module with $G$-action closely related to the irreducible representation $\rho$.  Under suitable additional hypotheses (e.g., $X,Y$ integral and Noetherian, $\pi$ flat), the sheaf $\F_{\rho}$ is locally free of rank equal to $(\text{dim }\rho)^2$. Conversely, to construct a cover given an appropriate collection of locally free $\Oh_X[G]$-modules $\{\F_{\rho}\}_{\rho}$, one must define a commutative, associative $\Oh_X[G]$-algebra structure on $\A=\oplus_{\rho} \F_{\rho}$.  One then obtains the $G$-cover $\pi:\RelSpec{X}{\A}\to X$.

The theory for abelian groups was analyzed by Pardini in \cite{pardini}.  In this case, the decomposition runs over the irreducible characters of $G$, and the $\Oh_X[G]$-submodule $\F_{\chi}$ is the invertible $\chi$-eigensheaf of $\pi_*\Oh_Y$, defined by the collections of sections on which the group acts as multiplication by the character $\chi$.  The algebra structure on $\pi_*\Oh_Y$ is determined by a compatible collection of morphisms $\{\F_{\chi}\otimes \F_{\chi'} \to \F_{\chi \chi'}\}_{\chi, \chi'}$, or equivalently, by a collection of global sections of $\{\F_{\chi}^{-1}\otimes\F_{\chi'}^{-1}\otimes\F_{\chi \chi'}\}_{\chi,\chi'}$. These sections are closely related to the branch divisor of the cover: given invertible sheaves $\{\F_{\chi}\}_{\chi}$, to construct a $G$-cover one may replace the explicit definition of the algebra structure with a specification of the branching data. As long as a  ``covering condition" is satisfied, one so builds a $G$-cover.

A key aspect of the abelian theory is that it allows one to understand the geometry of a covering scheme in terms of the geometry of the (usually simpler) base scheme. For example, if $X$ is a surface, one can use geometrically interesting configurations of curves in $X$ to construct new surfaces whose intrinsic geometry reflects the geometry of the configuration. This strategy has proven remarkably fruitful.  For example, a standard result in the theory of complex surfaces is the Bogomolov-Miyaoka-Yau inequality, which states that $\frac{K_X^2}{\chi ({\mathcal O}_X)}\leq 9$ for any smooth, complex surface $X$ of general type (\cite{bogomolov},\cite{miyaoka},\cite{yau}).
This inequality is sharp, and Hirzebruch produced examples of equality by constructing abelian covers of $\mathbb{P}^2$ branched over ``extreme" configurations of lines \cite{hirzebruch}. The inequality is known to fail in positive characteristic \cite{lang}, and infinite families of counterexamples can be produced using abelian covers branched over configurations of lines occuring only in positive characteristic \cite{easton}. A second example of this strategy can be found in \cite{eastonvakil}, in which a similar construction with abelian covers was employed to prove a higher-dimensional analogue of Belyi's theorem.

The situation in the nonabelian case is more complicated, perhaps even intractable in general, but the permutation group $S_3$ is within reach.  The aim here is to create structures parallel to those found in the abelian case, to the extent possible, with the hope of eventually extending the applications of the abelian theory to this nonabelian setting.  To realize this goal, we first proceed with a local analysis of $S_3$-covers, which we subsequently use as a guide in studying the global situation.  The precise statement of the main result (Theorem \ref{maintheorem}) requires some preliminary definitions, but can be roughly stated as follows:

\tpoint{Main Theorem}{\em Let $X$ be an integral, Noetherian scheme over an integral domain $R$ in which $6$ is invertible.  Then the collection of all flat $S_3$-covers of $X$ is parameterized by the following data:
\begin{itemize}
\item[(i)] an invertible sheaf $\mathcal{L}$ on $X$, on which $S_3$ acts via the sign character;
\item[(ii)] a locally free $\Oh_X[S_3]$-module $\E$, on which $S_3$ acts through its two-dimensional representation;
\item[(iii)] a module, $\text{Build}_X(\mathcal{L},\E)$, parameterizing the ``building data" which define commutative, associative algebra structures on $\A:=\Oh_X\oplus \mathcal{L}\oplus \E$ compatible with the given $S_3$-actions; i.e., the data required to construct a cover of the form $\pi:\RelSpec{X}{\A}\to X$.
\end{itemize}}

\noindent {\bf Acknowledgments.} The current work grew out of \cite{thesis}, simplifying the local analysis contained in the latter and extending the analysis to the global situation.  I am indebted to R. Vakil for his constant guidance and advice, both on the former project and in the present one.

\section{Preliminary Analysis}
\label{Prelim}

Fix a presentation
\[
S_3=\langle \sigma, \tau\; | \; \sigma^3=\tau^2=e, \quad \tau\sigma = \sigma^2\tau \rangle.
\]
Let $R$ be an integral domain in which $6$ is invertible.  The group ring $R[S_3]$ can be decomposed (as a free $R$-module) as $R[S_3]=C_1 \oplus C_2 \oplus C_3$, where
\begin{align*}
C_1&=\text{span}_R\{e+\sigma+\sigma^2+\tau+\sigma\tau+\sigma^2\tau\}=\{v\in R[S_3]\; | \; g\cdot v = v, \; \forall g\in S_3\},\\
C_2&=\text{span}_R\{e+\sigma+\sigma^2-\tau-\sigma\tau-\sigma^2\tau\}=\{v\in R[S_3]\; | \; g\cdot v = \text{sgn}(g)v, \; \forall g\in S_3\},
\end{align*}
and where a basis for $C_3$ is given by the four vectors
\begin{align*}
u_{11}&= -e+\sigma+\tau-\sigma^2\tau\\
u_{12}&=-\sigma+\sigma^2-\tau+\sigma\tau\\
u_{21}&= -e+\sigma^2+\tau-\sigma\tau\\
u_{22}&= e-\sigma+\sigma\tau-\sigma^2\tau.
\end{align*}
Under this decomposition of $R[S_3]$, we have $e=e_1+e_2+e_3$, where
\begin{align*}
e_1&=\frac{1}{6}(e+\sigma+\sigma^2+\tau+\sigma\tau+\sigma^2\tau),\\
e_2&=\frac{1}{6}(e+\sigma+\sigma^2-\tau-\sigma\tau-\sigma^2\tau),\\
e_3&= \frac{1}{3}(2e-\sigma-\sigma^2).
\end{align*}
Note that each $e_i$ is in the center of $R[S_3]$ and satisfies $e_ie_j=\delta_{ij}e_i$.  We also have a (non-equivariant) decomposition $C_3=C_{31}\oplus C_{32}$, where $C_{3i}=\text{span}_R\{u_{i1},u_{i2}\}$.  Under this decomposition, $e_3$ decomposes as $e_3=e_{31}+e_{32}$, where
\begin{align*}
e_{31}&=\frac{1}{3}(e-\sigma+\sigma\tau-\sigma^2\tau),\\
e_{32}&= \frac{1}{3}(e-\sigma^2-\sigma\tau +\sigma^2\tau),
\end{align*}
which also satisfy $e_{3i}e_{3j}=\delta_{ij}e_{3i}$ (but are not central).

Suppose $X$ is a scheme over $R$, and $\F$ is an $\Oh_X[S_3]$-module, with action explicitly given by a group homomorphism
\[
\mu: S_3\to \text{Aut}_{\Oh_X}(\F).
\]
Extend this morphism $R$-linearly to a ring homomorphism
\[
\mu: R[S_3]\to \text{End}_{\Oh_X}(\F).
\]
It is immediate that, given any $f\in \text{End}_{\Oh_X}(\F)$ satisfying $f\circ f = f$, the image presheaf of $f$ is automatically a sheaf.  Let $\F_i\subset \F$ denote the image (pre)sheaf of $\mu(e_i)$. Note that $\F_i$ is characterized by the property that $\mu(e_j)$ acts as $\delta_{ij}\text{Id}_{\F_i}$.
(In particular, $\F_1$ is the subsheaf of invariant sections of $\F$, and $\F_2$ is the subsheaf of sections on which $S_3$ acts as multiplication by the sign character.) Since $\mu(e)=\text{Id}_{\F}$, and since $e$ decomposes as $e=e_1+e_2+e_3$ with elements $e_i$ satisfying the aforementioned properties, the morphisms
$\mu(e_i)$ induce an $S_3$-equivariant direct sum decomposition
\[
\F=\F_1\oplus \F_2 \oplus \F_3.
\]
We can further decompose the sheaf $\F_3$, using the decomposition $e_3=e_{31}+e_{32}$.  We obtain an $\Oh_X$-module direct sum decomposition
\[
\F_3=\F_{31}\oplus \F_{32},
\]
where $\F_{3i}$ is the image (pre)sheaf associated to $\mu(e_{3i})|_{\F_3}$.  Note that these subsheaves are not invariant under the group action.
Indeed, observe that
\[
\tau e_{31}=\frac{1}{3}(\tau -\sigma^2\tau +\sigma^2-\sigma)=e_{32}\tau,
\]
so $\mu(\tau)|_{\F_3}$ is an automorphism interchanging the summands. We can therefore write $\F_3=\F_{31}\oplus \tau \F_{31}$.

Suppose $\pi:Y\to X$ is a flat $S_3$-cover of $R$-schemes.  Then $\pi_*\Oh_Y$ is a locally free $\Oh_X[S_3]$-module of rank $6$, with $(\pi_*\Oh_Y)^{S_3}=\Oh_X$. By the above, we have an induced $\Oh_X[S_3]$-module direct sum decomposition $\pi_*\Oh_Y=\Oh_X\oplus \F_2\oplus \F_3$ (as well as an $\Oh_X$-module decomposition $\F_3=\F_{31}\oplus \tau \F_{31}$).  If we assume $X, Y$ are integral, Noetherian $R$-schemes, then it follows that $\F_2$ and $\F_{31}$ are also locally free, of ranks $1$ and $2$, respectively.

Since a finite morphism is affine, it follows that $Y=\RelSpec{X}{\left(\Oh_X\oplus \F_2\oplus \F_3\right)}$.  Thus, to construct a flat $S_3$-cover of a given scheme $X$, one needs the following data:
\begin{itemize}
\item[(i)] An invertible sheaf $\mathcal{L}$, on which $S_3$ acts via the sign character;
\item[(ii)] A locally free $\Oh_X$-module $\E$ of rank $4$, together with an $S_3$-action such that (under the induced $R[S_3]$-action), $e_i$ acts as $\delta_{3i}\text{Id}_{\E}$; and
\item[(iii)] A commutative, associative $\Oh_X$-algebra structure on $\A:=\Oh_X\oplus \mathcal{L} \oplus \E$ compatible with the given $S_3$-action.
\end{itemize}

We aim to precisely describe data (iii), given data (i) and (ii).  To that end, suppose we are given such $\Oh_X$-modules $\mathcal{L}$ and $\E$.

\tpoint{Lemma} {\em A commutative $\Oh_X$-algebra structure on $\A$ compatible with the given $S_3$-actions is defined by a triple of $\Oh_X$-module morphisms
\[ \alpha: S^2\mathcal{L}\to \Oh_X, \quad \beta: \mathcal{L}\otimes_{\Oh_X}\E\to \E, \quad \gamma: S^2\E\to \A.\]}

\begin{proof}
A priori, an algebra structure is given by an $\Oh_X$-module morphism $\A\otimes_{\Oh_X}\A\to \A$, which is equivalent to an $\Oh_X$-module morphism
\begin{align*}
(\Oh_X&\otimes_{\Oh_X}\Oh_X)\oplus (\Oh_X\otimes_{\Oh_X}\mathcal{L})\oplus (\mathcal{L}\otimes_{\Oh_X}\Oh_X)\oplus(\Oh_X\otimes_{\Oh_X}\E)\oplus(\E\otimes_{\Oh_X}\Oh_X)\\
&\oplus(\mathcal{L}\otimes_{\Oh_X}\mathcal{L})\oplus(\mathcal{L}\otimes_{\Oh_X}\E)\oplus(\E\otimes_{\Oh_X}\mathcal{L})\oplus(\E\otimes_{\Oh_X}\E)\to A.
\end{align*}
The first coordinate of this morphism must be given by the algebra structure on $\Oh_X$,
and the following four coordinates must be given by the left and right $\Oh_X$-module structures on $\mathcal{L}$ and $\E$, respectively.  Commutativity requires the sixth and ninth coordinates to factor through the canonical morphisms to the corresponding symmetric products:
\[
\xymatrix{
\mathcal{L}\otimes_{\Oh_X}\mathcal{L} \ar[r] \ar[d]_-{\text{can}} & \A \\
S^2\mathcal{L} \ar[ur]
}, \quad
\xymatrix{
\E\otimes_{\Oh_X}\E \ar[r] \ar[d]_-{\text{can}} & \A \\
S^2\E \ar[ur]_{\gamma}
}.
\]
Similarly, the seventh and eighth coordinates must agree after the canonical braid isomorphism:
\[
\xymatrix{
\mathcal{L}\otimes_{\Oh_X}\E \ar[d]_{\cong}^-{\text{can}} \ar[r] & \A \\
\E\otimes_{\Oh_X}\mathcal{L} \ar[ur]
}.
\]

Compatibility with the $S_3$-action implies that $S_3$ acts as $(\text{sgn})^2=\text{Id}$ on the image of $S^2\mathcal{L}\to \A$, and so the morphism must factor
through the subsheaf $\Oh_X\subset \A$ of invariant sections:
\[
\xymatrix{
S^2\mathcal{L} \ar[dr]_{\alpha} \ar[r] & \A \\
& \Oh_X \ar@{^{(}->}[u]
}.
\]
Similarly, $\mathcal{L}\otimes_{\Oh_X}\E \to \A$ must factor through $\E\subset \A$:
\[
\xymatrix{
\mathcal{L}\otimes_{\Oh_X}\E \ar[dr]_{\beta} \ar[r] & \A \\
& \E \ar@{^{(}->}[u]
}.
\]
(Suppose $U\subset X$ is an affine open and $t\in \mathcal{L}(U), x\in \E(U)$.  Then
\begin{align*}
\mu(e_3) (t\otimes x) &= \frac{1}{3}(2\mu(e)(t\otimes x)-\mu(\sigma) (t\otimes x)-\mu(\sigma^2) (t\otimes x))\\
&=\frac{1}{3}(2(t\otimes x)-(\mu(\sigma) t)\otimes (\mu(\sigma) x)-(\mu(\sigma^2) t)\otimes (\mu(\sigma^2) x))\\
&= \frac{1}{3}(t\otimes(2x)-t\otimes(\mu(\sigma) x)-t\otimes (\mu(\sigma^2) x))\\
&= t\otimes (\mu(e_3) x)\\
&= t\otimes x.
\end{align*}
Similarly, one sees $\mu(e_1)(t\otimes x)=t\otimes (\mu(e_2)x)=0$ and $\mu(e_2)(t\otimes x)=t\otimes (\mu(e_1)x)=0$.)
\end{proof}

Of course, not every such triple of morphisms will satisfy all of the necessary compatibilities with the $S_3$-action, nor will it necessarily
define an associative algebra structure.  To precisely identify which triples do satisfy these conditions, we first analyze the local situation.

\section{Local Analysis}
\label{localanalysis}

Continuing the notation above, let $\E=\E'\oplus \tau \E'$ be the $\Oh_X$-module decomposition induced by $e_3=e_{31}+e_{32}$. Let $U\subset X$ be any affine open such that $\mathcal{L}(U),\E'(U)$ are free $\Oh_X(U)$-modules of ranks $1$, $2$, respectively. For notational simplicity, write $L=\mathcal{L}(U), E=\E(U), E'=\E'(U), B=\Oh_X(U)$, and $A=\A(U)$. Let $\{v_1, v_2\}$ be any basis for $E'$ (and so $\{v_1,v_2,\tau v_1, \tau v_2\}$ is a basis for $E$).

\tpoint{Lemma}\label{actionlemma}{\em
The $S_3$-action on $E$ satisfies:
\begin{enumerate}
\item[(i)] $\sigma v_i = -\tau v_i$
\item[(ii)] $\sigma^2 v_i = -v_i+\tau v_i$
\item[(iii)] $\sigma\tau v_i = v_i - \tau v_i$
\item[(iv)] $\sigma^2\tau v_i = -v_i$.
\end{enumerate}
}

\begin{proof}
By direct calculation, one sees
\begin{align*}
\sigma e_{31}&=\frac{1}{3}(\sigma-\sigma^2+\sigma^2\tau-\tau)=-\tau e_{31},\\
\sigma^2 e_{31}&=\frac{1}{3}(\sigma^2-e+\tau-\sigma\tau)=(-e+\tau)e_{31}.
\end{align*}
It follows that $\mu(\sigma)|_{E'}=\mu(-\tau)|_{E'}$ and $\mu(\sigma^2)|_{E'}=\mu(-e+\tau)|_{E'}$, which prove (i) and (ii).
Properties (iii) and (iv) then immediately follow from the relations $\sigma\tau=\tau\sigma^2$ and $\sigma^2\tau=\tau\sigma$.
\end{proof}

Let $t\in L$ be any generator.  To understand the algebra structure, we need to analyze (the images of) the following products (with tensor symbols suppressed):
\begin{itemize}
\item $t^2\in B$
\item $tv_i, t\cdot\tau v_i \in E$, for $i=1, 2$
\item $v_iv_j, v_i \cdot\tau v_j, \tau v_i \cdot\tau v_j\in A$, for $i,j=1,2$
\end{itemize}

\tpoint{Lemma}\label{preassocstruct} {\em With respect to the basis $\{1, t, v_1, v_2, \tau v_1, \tau v_2\}$ for $A$, the commutative $B[S_3]$-algebra structures on $A$
are precisely those of the form
\begin{align*}
t^2&= a\\
tv_1&= b_1v_1+b_2v_2-2b_1\tau v_1 - 2b_2\tau v_2\\
tv_2&= c_1v_1+c_2v_2-2c_1\tau v_1 - 2c_2\tau v_2\\
t\cdot \tau v_1&= 2b_1v_1+2b_2v_2-b_1\tau v_1-b_2\tau v_2\\
t\cdot \tau v_2&= 2c_1v_1+2c_2v_2-c_1\tau v_1-c_2\tau v_2\\
v_1^2&= d_1+d_3v_1+d_4v_2-2d_3\tau v_1 - 2d_4\tau v_2\\
v_1v_2&= e_1+e_3v_1+e_4v_2-2e_3\tau v_1 - 2e_4\tau v_2\\
v_2^2&= f_1+f_3v_1+f_4v_2-2f_3\tau v_1-2f_4\tau v_2\\
v_1\cdot \tau v_1&= \frac{1}{2}d_1-d_3v_1-d_4v_2-d_3\tau v_1-d_4\tau v_2\\
v_1\cdot \tau v_2&= \frac{1}{2}e_1+h_2t+h_3v_1+h_4v_2-e_3\tau v_1 -e_4\tau v_2\\
v_2\cdot \tau v_1&= \frac{1}{2}e_1-h_2t-e_3v_1-e_4v_2+h_3\tau v_1+h_4\tau v_2\\
v_2\cdot \tau v_2&= \frac{1}{2}f_1-f_3v_1-f_4v_2-f_3\tau v_1-f_4\tau v_2\\
(\tau v_1)^2&= d_1-2d_3v_1-2d_4v_2+d_3\tau v_1+d_4\tau v_2\\
(\tau v_1)(\tau v_2)&= e_1-2e_3v_1-2e_4v_2+e_3\tau v_1 +e_4\tau v_2\\
(\tau v_2)^2&= f_1-2f_3v_1-2f_4v_2+f_3\tau v_1 + f_4\tau v_2,
\end{align*}
for some $a, b_i, c_i, d_i, e_i, f_i, h_i \in B$.
}

\begin{proof} Note that compatibility with the action of $\tau$ requires
\[
\tau(tv_i)=(\tau t)(\tau v_i)=-t \cdot\tau v_i,
\]
and so $t \cdot\tau v_i$ is determined by $t v_i$.  Similarly, $\tau v_i \cdot\tau v_j$ is determined by $v_i v_j$,
and $v_2 \cdot \tau v_1$ is determined by $v_1 \cdot \tau v_2$.  So, suppose
\begin{align*}
t^2 &= a\\
tv_1 &= b_1v_1 + b_2v_2 + b_3\tau v_1 + b_4\tau v_2\\
tv_2 &= c_1v_1 + c_2v_2 + c_3\tau v_1 + c_4\tau v_2\\
v_1^2 &= d_1 + d_2t + d_3v_1 + d_4v_2 + d_5\tau v_1 + d_6\tau v_2\\
v_1v_2 &= e_1 + e_2t + e_3v_1 + e_4v_2 + e_5\tau v_1 + e_6\tau v_2\\
v_2^2 &= f_1 + f_2t + f_3v_1 + f_4v_2 + f_5\tau v_1 + f_6\tau v_2\\
v_1\cdot \tau v_1 &= g_1 + g_2t + g_3v_1 + g_4v_2 + g_5\tau v_1 + g_6\tau v_2\\
v_1\cdot \tau v_2 &= h_1 + h_2t + h_3v_1 + h_4v_2 + h_5\tau v_1 + h_6\tau v_2\\
v_2\cdot \tau v_2 &= i_1 + i_2t + i_3v_1 + i_4v_2 + i_5\tau v_1 + i_6\tau v_2,
\end{align*}
for some $a, b_j, c_j, d_j, e_j, f_j, g_j, h_j, i_j\in B$.
By our previous remark, we then must have
\begin{align*}
t\cdot \tau v_1 &= -b_3v_1-b_4v_2-b_1\tau v_1 - b_2\tau v_2\\
t\cdot \tau v_2 &= -c_3v_1-c_4v_2 -c_1\tau v_1 - c_2\tau v_2\\
(\tau v_1)^2 &= d_1 - d_2t +d_5v_1 +d_6v_2 + d_3\tau v_1 + d_4\tau v_2\\
(\tau v_1)(\tau v_2) &= e_1 - e_2t +e_5v_1 +e_6v_2 + e_3\tau v_1 + e_4\tau v_2\\
(\tau v_2)^2 &= f_1 - f_2t +f_5v_1 +f_6v_2 + f_3\tau v_1 + f_4\tau v_2\\
v_2\cdot \tau v_1 &= h_1 - h_2t +h_5v_1 +h_6v_2 + h_3\tau v_1 + h_4\tau v_2.
\end{align*}
Such an algebra structure is now ensured to be compatible with the action of $\tau$, and so compatibility with the full $S_3$-action will follow from
compatibility with the action of $\sigma$.  We now impose this compatibility on each product.  Using the equations above,
together with Lemma \ref{actionlemma}, we compute
\begin{align*}
\sigma(tv_1)&= b_3v_1+b_4v_2+(-b_1-b_3)\tau v_1 + (-b_2-b_4)\tau v_2,\\
\intertext{and}
(\sigma t)(\sigma v_1)&= t\cdot \sigma v_1 = t\cdot (-\tau v_1)= -(t\cdot \tau v_1)\\
&= b_3v_1 + b_4v_2 + b_1\tau v_1 + b_2\tau v_2.
\end{align*}
So, compatibility with $\sigma$ requires $b_3=-2b_1$ and $b_4=-2b_2$.  The corresponding computation for the relation $\sigma(tv_2)=(\sigma t)(\sigma v_2)$ requires $c_3=-2c_1$ and $c_4 = -2c_2$.

Similarly, we compute
\begin{align*}
\sigma(v_1^2)&= d_1 + d_2t + d_5v_1 + d_6v_2 + (-d_3-d_5)\tau v_1 + (-d_4-d_6)\tau v_2,\\
\intertext{and}
(\sigma v_1)^2&= (-\tau v_1)^2 = (\tau v_1)^2\\
&= d_1 - d_2t + d_5v_1 + d_6v_2 + d_3\tau v_1 + d_4\tau v_2.
\end{align*}
So, we must have $d_2=0$, $d_5=-2d_3$, and $d_6=-2d_4$.  Similarly, the relation $\sigma(v_1v_2)=\sigma(v_1)\sigma(v_2)$ requires
$e_2=0$, $e_5=-2e_3$, and $e_6=-2e_4$, while the relation $\sigma(v_2)^2=(\sigma v_2)^2$ requires $f_2=0$, $f_5=-2f_3$ and $f_6=-2f_4$.

Lastly, we compute
\begin{align*}
\sigma(v_1\cdot \tau v_1)&= g_1 + g_2t + g_5v_1 + g_6v_2 + (-g_3-g_5)\tau v_1 + (-g_4-g_6)\tau v_2,\\
\intertext{and}
(\sigma v_1)(\sigma \tau v_1)&= (-\tau v_1)(v_1-\tau v_1)=-(v_1\cdot \tau v_1)+(\tau v_1)^2\\
&= (-g_1+d_1) + (-g_2-d_2)t + (-g_3+d_5)v_1 + (-g_4+d_6)v_2 \\
& \quad + (-g_5+d_3)\tau v_1 + (-g_6+d_4)\tau v_2.
\end{align*}
From this (and the previously obtained relations), we deduce $g_1=\frac{1}{2}d_1$, $g_2=0$, $g_3=g_5=-d_3$, and $g_4=g_6=-d_4$.
Similarly, the relation
$\sigma(v_1\cdot \tau v_2)=(\sigma v_1)(\sigma \tau v_2)$ requires $h_1=\frac{1}{2}e_1$, $h_5=-e_3$, and $h_6=-e_4$, while the relation
$\sigma(v_2\cdot \tau v_2)=(\sigma v_2)(\sigma\tau v_2)$ requires $i_1=\frac{1}{2}f_1$, $i_2=0$, $i_3=i_5=-f_3$ and $i_4=i_6=-f_4$.
\end{proof}

Notice we specifically omitted any mention of associativity.  Indeed, the associativity relations impose many additional conditions on the system.

\tpoint{Proposition}\label{localstructure} {\em With respect to the basis $\{1, t, v_1, v_2, \tau v_1, \tau v_2\}$ for $A$, every commutative, associative, integral
$B[S_3]$-algebra structure on $A$ is of the form
\begin{align*}
t^2&= -3b_1^2-3b_2c_1\\
tv_1&= b_1v_1+b_2v_2-2b_1 \tau v_1 -2b_2\tau v_2\\
tv_2&= c_1v_1-b_1v_2-2c_1\tau v_1 +2b_1\tau v_2\\
t\cdot \tau v_1&= 2b_1v_1+2b_2v_2-b_1\tau v_1 -b_2\tau v_2\\
t\cdot \tau v_2&= 2c_1v_1-2b_1v_2-c_1\tau v_1 +b_1\tau v_2\\
v_1^2&= 6(d_3^2-d_4f_4)+d_3v_1+d_4v_2-2d_3\tau v_1 -2d_4\tau v_2\\
v_1v_2&=3(d_4f_3-d_3f_4)-f_4v_1-d_3v_2+2f_4\tau v_1 + 2d_3\tau v_2\\
v_2^2&= 6(f_4^2-d_3f_3)+f_3v_1+f_4v_2-2f_3\tau v_1 -2f_4\tau v_2\\
v_1\cdot \tau v_1 &= 3(d_3^2-d_4f_4)-d_3v_1-d_4v_2-d_3\tau v_1-d_4\tau v_2\\
v_1\cdot \tau v_2&= \frac{3}{2}(d_4f_3-d_3f_4)+h_2t+f_4v_1+d_3v_2+f_4\tau v_1 +d_3\tau v_2\\
v_2\cdot \tau v_1&= \frac{3}{2}(d_4f_3-d_3f_4)-h_2t+f_4v_1+d_3v_2+f_4\tau v_1 + d_3\tau v_2\\
v_2\cdot \tau v_2&= 3(f_4^2-d_3f_3)-f_3v_1-f_4v_2-f_3\tau v_1 -f_4\tau v_2\\
(\tau v_1)^2&= 6(d_3^2-d_4f_4)-2d_3v_1-2d_4v_2+d_3\tau v_1+d_4\tau v_2\\
(\tau v_1)(\tau v_2)&= 3(d_4f_3-d_3f_4)+2f_4v_1+2d_3v_2-f_4\tau v_1-d_3\tau v_2\\
(\tau v_2)^2&= 6(f_4^2-d_3f_3)-2f_3v_1-2f_4v_2+f_3\tau v_1 + f_4\tau v_2,
\end{align*}
for some $b_1, b_2, c_1, d_3, d_4, f_3, f_4, h_2 \in B$ satisfying
\begin{enumerate}
\item [(i)] $2b_1d_3-b_2f_4+c_1d_4=0$;
\item [(ii)] $2b_1f_4-b_2f_3+c_1d_3=0$; and
\item [(iii)] $(b_1^2+b_2c_1)h_2=\frac{3}{2}(b_1(d_4f_3-d_3f_4)+b_2(f_4^2-d_3f_3)+c_1(d_4f_4-d_3^2))$.
\end{enumerate}

Conversely, any multiplicative structure on $A$ of the above form defines a commutative, associative (but possibly non-integral) $B[S_3]$-algebra structure on $A$.
}

\begin{proof}
The proof consists of systematically imposing the third-order associativity conditions.  Using Lemma \ref{preassocstruct}, we compute
\begin{align*}
(t^2)v_1&= av_1\\
t(tv_1)&= (-3b_1^2-3b_2c_1)v_1+(-3b_1b_2-3b_2c_2)v_2.
\end{align*}
Equating coefficients gives
\begin{align}
a&= -3b_1^2-3b_2c_1\\
0&= b_2(b_1+c_2).
\end{align}
Similarly, we compute
\begin{align*}
(t^2)v_2&= av_2\\
t(tv_2)&= (-3b_1c_1-3c_1c_2)v_1+(-3b_2c_1-3c_2^2)v_2,
\end{align*}
and so must have
\begin{align}
0&= c_1(b_1+c_2)\\
a&=-3b_2c_1-3c_2^2.
\end{align}
Note that the relations $(t^2)\tau v_1=t(t\tau v_1)$ and $(t^2)\tau v_2=t(t\tau v_2)$ immediately follow from the above relations and the compatibility with $\tau$.
Indeed, we have $(t^2)\tau v_1=(-t)^2\tau v_1=(\tau t)^2\tau v_1=\tau( (t^2)v_1)=\tau (t(tv_1))=-t(-t\cdot \tau v_1)=t(t\tau v_1)$, and similarly for $(t^2)\tau v_2$.

We next compute
\begin{align*}
t(v_1^2)&= d_1t+(-3b_1d_3-3c_1d_4)v_1+(-3b_2d_3-3c_2d_4)v_2,\\
(tv_1)v_1&= (-2b_2h_2)t+(3b_1d_3+b_2e_3-2b_2h_3)v_1+(3b_1d_4+b_2e_4-2b_2h_4)v_2,
\end{align*}
which implies
\begin{align}
d_1&=-2b_2h_2\\
6b_1d_3+b_2e_3-2b_2h_3+3c_1d_4&=0\\
3b_1d_4+3b_2d_3+b_2e_4-2b_2h_4+3c_2d_4&=0.
\end{align}

We claim integrality of the algebra structure requires $c_2=-b_1$.
Indeed, suppose $c_2\neq -b_1$.  Then equations (2) and (3) imply $b_2=c_1=0$.  Equations (1) and (4) then become $-3c_2^2=a=-3b_1^2$.
By hypothesis, the algebra is integral, so $t^2=a$ is nonzero, and hence it follows that $b_1$ and $c_2$ are both nonzero.  Since their squares are equal and $c_2\neq -b_1$, we must have $c_2=b_1$.  But then equations $(5)-(7)$ imply $d_1=d_3=d_4=0$, and hence $v_1^2=0$, which violates integrality.

So, we must have $c_2=-b_1$, and equations (1)-(7) now reduce to
\begin{align}
a&= -3b_1^2-3b_2c_1\\
d_1&= -2b_2h_2\\
6b_1d_3+b_2e_3-2b_2h_3+3c_1d_4&=0\\
3b_2d_3+b_2e_4-2b_2h_4&=0.
\end{align}

We now continue computing the associativity relations:
\begin{align*}
t(v_1v_2)&= e_1t+(-3b_1e_3-3c_1e_4)v_1+(-3b_2e_3+3b_1e_4)v_2\\
(tv_1)v_2&= (2b_1h_2)t+(3b_1e_3+3b_2f_3)v_1+(3b_1e_4+3b_2f_4)v_2\\
&\quad +(-2b_1e_3-2b_1h_3)\tau v_1 +(-2b_1e_4-2b_1h_4)\tau v_2.
\end{align*}
It follows that
\begin{align}
e_1&=2b_1h_2\\
2b_1e_3+b_2f_3+c_1e_4&=0\\
b_2(e_3+f_4)&=0\\
b_1(e_3+h_3)&=0\\
b_1(e_4+h_4)&=0.
\end{align}

We next compute:
\begin{align*}
t(v_2^2)&= f_1t+(-3b_1f_3-3c_1f_4)v_1+(-3b_2f_3+3b_1f_4)v_2\\
(tv_2)v_2&= (2c_1h_2)t+(3c_1e_3-3b_1f_3)v_1+(3c_1e_4-3b_1f_4)v_2\\
&\quad +(-2c_1e_3-2c_1h_3)\tau v_1 +(-2c_1e_4-2c_1h_4)\tau v_2,
\end{align*}
which implies
\begin{align}
f_1&= 2c_1h_2\\
c_1(e_3+f_4)&= 0\\
2b_1f_4-b_2f_3-c_1e_4&=0\\
c_1(e_3+h_3)&=0\\
c_1(e_4+h_4)&=0.
\end{align}

Now, again by integrality, we cannot have both $b_1$ and $b_2$ both zero (else $tv_1=0$), nor $b_1$ and $c_1$ both zero (else $tv_2=0$).
In particular, equations (15) and (20) together imply $h_3=-e_3$.  Similarly, equations (16) and (21) imply $h_4=-e_3$.
Also observe that summing equations (13) and (19) yields the equation $b_1(e_3+f_4)$, which together with equation (14) (or (18)) implies $e_3=-f_4$.

We once again continue with the associativity relations.  A calculation reveals the relations
\[
t(v_1\tau v_1)=(tv_1)\tau v_1, \quad t(v_1\tau v_2)=(tv_1)\tau v_2,\quad t(v_2\tau v_2)=(tv_2)\tau v_2
\]
all now hold (under the conditions summarized above).  Compatibility with the action of $\tau$ then implies that the relations
\[
t(\tau v_1)^2=(t\cdot \tau v_1)\tau v_1,\quad t(\tau v_1\cdot \tau v_2)=(t\cdot \tau v_1)\tau v_2,\quad t(\tau v_2)^2=(t\cdot \tau v_2)\tau v_2, \quad t(v_2\cdot \tau v_1)
=(t v_2)\tau v_1
\]
automatically follow from the previous associativity relations.  So, we next compute
\begin{align*}
(v_1^2)v_2&= (2d_3h_2)t+(-3d_3f_4+3d_4f_3)v_1+(-2b_2h_2+3d_3e_4+3d_4f_4)v_2\\
v_1(v_1v_2)&= (-2e_4h_2)t+(2b_1h_2-3d_3f_4-3e_4f_4)v_1+(-3d_4f_4+3e_4^2)v_2,
\end{align*}
which implies
\begin{align}
h_2(d_3+e_4)&=0\\
2b_1h_2-3d_4f_3-3e_4f_4&=0\\
2b_2h_2-3d_3e_4-6d_4f_4+3e_4^2&=0.
\end{align}

Continuing, we compute
\begin{align*}
v_1(v_2^2)&= (-2f_4h_2)t+(2c_1h_2+3d_3f_3-3f_4^2)v_1+(d_4f_2+3e_4f_4)v_1\\
(v_1v_2)v_2&= (-2f_4h_2)t+(3f_4^2+3e_4f_3)v_1+(2b_1h_2)v_2,
\end{align*}
which implies
\begin{align}
2c_1h_2+3d_3f_3-3e_4f_3-6f_4^2&=0.
\end{align}

Lastly, we compute
\begin{align*}
v_1^2(\tau v_1)&= (3b_2d_3h_2-3b_1d_4h_2)+(-d_4h_2)t+(3d_3^2-3d_4f_4)v_1+(3d_3d_4+3d_4e_4)v_2\\
&\quad +(-2b_2h_2-3d_3^2+3d_4f_4)\tau v_1+(-3d_3d_4-3d_4e_4)\tau v_2\\
v_1(v_1\tau v_1)&= (3b_2d_3h_2-3b_1d_4h_2)+(-d_4h_2)t+(-b_2h_2)v_1\\
&\quad +(3d_3^2-3d_4f_4)\tau v_1+(3d_3d_4+3d_4e_4)\tau v_2,
\end{align*}
which implies
\begin{align}
b_2h_2+3d_3^2-3d_4f_4&=0\\
d_4(d_3+e_4)&=0.
\end{align}

Observe that if $e_4\neq -d_3$, then equations (11), (22), and (27) imply $b_2=h_2=d_4=0$, which together with equation (26) imply $d_3=0$.
But then $v_1^2=0$, which violates integrality.  So, we must have $e_4=-d_3$.

At this point, we've reached the statement of Proposition \ref{localstructure}.  Indeed, equations $(23)-(25)$ (together with $(5),(12),(17)$) can all be immediately solved to give the constant coefficients
of the terms $v_1^2, \ldots , (\tau v_2)^2$, as well as combined to give relation (iii) of the proposition.  A calculation verifies all remaining associativity relations
are now satisfied.
\end{proof}

\tpoint{Corollary}{\em Continuing the notation from Proposition \ref{localstructure}, the ramification locus of $\pi:\Spec{A}\to \Spec{B}$ is the zero locus of the ideal generated by all $5$ x $5$ minors of the matrix
\[
\left[\begin{array}{ccccc}
2t & 0 & 0 & 0 & 0\\
v_1 & t-b_1 & -b_2 & 2b_1 & 2b_2\\
v_2 & -c_1 & t+b_1 & 2c_1 & -2b_1\\
\tau v_1 & -2b_1 & -2b_2 & t+b_1 & b_2\\
\tau v_2 & -2c_1 & 2b_1 & c_1 & t-b_1\\
0 & 2v_1-d_3 & -d_4 & 2d_3 & 2d_4\\
0 & v_2+f_4 & v_1+d_3 & -2f_4 & -2d_3\\
0 & -f_3 & 2v_2-f_4 & 2f_3 & 2f_4\\
0 & \tau v_1+d_3 & d_4 & v_1+d_3 & d_4\\
-h_2 & \tau v_2-f_4 & -d_3 & -f_4 & v_1-d_3\\
h_2 & -f_4 & \tau v_1-d_3 & v_2-f_4 & -d_3\\
0 & f_3 & \tau v_2+f_4 & f_3 & v_2+f_4\\
0 & 2d_3 & 2d_4 & 2\tau v_1 -d_3 & -d_4\\
0 & -2f_4 & -2d_3 & \tau v_2+f_4 & \tau v_1+d_3\\
0 & 2f_3 & 2f_4 & -f_3 & 2\tau v_2-f_4
\end{array}\right],
\]
under the identification of $A$ with the quotient of $B[1,t,v_1,v_2,\tau v_1, \tau v_2]$ by the ideal of relations generated by the equations of Proposition \ref{localstructure}.}

\tpoint{Corollary}\label{localtoglobal}{\em
The multiplication in $A$ is determined by a triple of morphisms
\[
\phi: L\otimes E'\to E' \quad \psi:S^2E'\to E' \quad \xi:E'\otimes \tau E'\to L.
\]
If $t$ is a generator for $L$ and $\{v_1, v_2\}$ is a basis for $E'$, then these morphisms are of the form
\begin{align*}
\phi(t\otimes v_1)&= av_1+bv_2\\
\phi(t\otimes v_2)&= cv_1-av_2\\
\psi(v_1^2)&= dv_1+ev_2\\
\psi(v_1v_2)&= -gv_1-dv_2\\
\psi(v_2^2)&= fv_1+gv_2\\
\xi(v_1\otimes \tau v_1)&= 0\\
\xi(v_1\otimes \tau v_2)&= ht\\
\xi(v_2\otimes \tau v_1)&= -ht\\
\xi(v_2\otimes \tau v_2)&= 0.
\end{align*}
for $a, b, c, d, e, f, g, h\in B$ satisfying the relations
    \begin{itemize}
    \item [(i)] $-bg+2ad+ce=0$;
    \item [(ii)] $-bf+2ag+cd=0$;
    \item [(iii)] $(a^2+bc)h=\frac{3}{2}(a(ef-dg)+b(g^2-df)+c(eg-d^2))$.
    \end{itemize}
}

\bpoint{Remark} The morphism $\psi$ is what Miranda calls a {\em triple cover homomorphism} in \cite{miranda}.  Such a homomorphism induces a triple cover $\RelSpec{X}{(\Oh_X\oplus \E')}\to X$.

\bpoint{Remark} There are certainly many solutions to the above system of constraints.  For instance, $a=b=c=d=g=1, e=-1, f=3, h=-6$ and $a=b=c=d=f=1,e=-2,g=0,h=-3$ are both possible solutions.

\section{Global Analysis}
\label{globalanalysis}

The previous section obtained a local description of flat $S_3$-covers.  We now use this local description to obtain a global description. In light of Corollary \ref{localtoglobal}, we now expect to characterize such covers by a submodule of
\[
\Hom(\mathcal{L}\otimes \E',\E')\oplus \Hom(S^2\E',\E')\oplus \Hom(\E'\otimes \tau \E',\mathcal{L}).
\]
We need a basis-free restatement of Corollary \ref{localtoglobal}.

\tpoint{Definition}
Let $M_1 \leq  \Hom(L\otimes E',E')$ denote the submodule consisting of elements $\phi$ of the form
\begin{align*}
\phi(t\otimes v_1)&= av_1+bv_2\\
\phi(t\otimes v_2)&= cv_1-av_2.
\end{align*}
Let $M_2 \leq \Hom(S^2E', E')$ denote the submodule consisting of elements $\psi$ of the form
\begin{align*}
\psi(v_1^2)&= dv_1+ev_2\\
\psi(v_1v_2)&= -gv_1-dv_2\\
\psi(v_2^2)&= fv_1+gv_2.
\end{align*}
Let $M_3 \leq \Hom(E'\otimes \tau E',L)$ denote the submodule consisting of elements $\xi$ of the form
\begin{align*}
\xi(v_1\otimes \tau v_1)&= 0\\
\xi(v_1\otimes \tau v_2)&= ht\\
\xi(v_2\otimes \tau v_1)&= -ht\\
\xi(v_2\otimes \tau v_2)&= 0.
\end{align*}

\tpoint{Lemma}{\em $M_1,M_2,$ and $M_3$ are well-defined.}

\begin{proof}
Let $s = at$ be another generator for $L$, and $\{w_1,w_2\}$ be another basis for $E'$, with change of basis matrix
\[
C=\left[\begin{array}{cc}\lambda_1 & \mu_1 \\ \lambda_2 &\mu_2\end{array}\right].
\]
A straightforward calculation then gives
\begin{align*}
\phi(s\otimes w_1)&= a'w_1+b'w_2\\
\phi(s\otimes w_2)&= c'w_1-a'w_2\\
\psi(w_1^2)&= d'w_1+e'w_2\\
\psi(w_1w_2)&= -g'w_1-d'w_2\\
\psi(w_2^2)&= f'w_1+g'w_2\\
\xi(w_1\otimes \tau w_1)&= 0\\
\xi(w_1\otimes \tau w_2)&= h's\\
\xi(w_2\otimes \tau w_1)&= -h's\\
\xi(w_2\otimes \tau w_2)&= 0,
\end{align*}
where
\begin{align*}
a'&= \frac{a}{\det(C)}(-\lambda_1\lambda_2b+\lambda_1\mu_2a+\lambda_2\mu_1a+\mu_1\mu_2c)\\
b'&= \frac{a}{\det(C)}(\lambda_1^2b-2\lambda_1\mu_1a-\mu_1^2c)\\
c'&= \frac{a}{\det(C)}(-\lambda_2^2b+2\lambda_2\mu_2a+\mu_2^2c)\\
d'&= \frac{1}{\det(C)}(-\lambda_1^2\lambda_2e+\lambda_1^2\mu_2d+2\lambda_1\lambda_2\mu_1d-2\lambda_1\mu_1\mu_2g-\lambda_2\mu_1^2g+\mu_1^2\mu_2f)\\
e'&= \frac{1}{\det(C)}(\lambda_1^3e-3\lambda_1^2\mu_1d+3\lambda_1\mu_1^2g-\mu_1^3f)\\
f'&= \frac{1}{\det(C)}(-\lambda_2^3e+3\lambda_2^2\mu_2d-3\lambda_2\mu_2^2g+\mu_2^3f)\\
g'&= \frac{1}{\det(C)}(\lambda_1\lambda_2^2e-2\lambda_1\lambda_2\mu_2d+\lambda_1\mu_2^2g-\lambda_2^2\mu_1d+2\lambda_2\mu_1\mu_2g-\mu_1\mu_2^2f)\\
h'&= \frac{\det(C)}{a}h.
\end{align*}
\end{proof}

\tpoint{Lemma}\label{basis-free-description}{\em There exist natural isomorphisms
\begin{itemize}
    \item[(i)] ${\bf F_1}:M_1 \tilde{\to} \Hom(L\otimes S^2E', \bigwedge^2E')$;
    \item[(ii)] ${\bf F_2}:M_2 \tilde{\to} \Hom(S^3E', \bigwedge^2E')$;
    \item[(iii)] ${\bf F_3}:M_3 \tilde{\to} \Hom(\bigwedge^2 E', L)$.
\end{itemize}}
\begin{proof}
The proof of (ii) is detailed in \cite[Prop. 3.3]{miranda}.  The proofs of (i) and (iii) are similar to that of (ii), and are given here. The method of proof will be used repeatedly.

We begin by proving (iii).  Assume $\xi \in M_3$ is of the form of Definition 5.2.  The induced morphism
\[
\xymatrix{
E'\otimes E' \ar[r]^-{1\otimes \tau} & E'\otimes \tau E' \ar[r]^-{\xi} & L
}
\]
maps
\begin{align*}
v_1\otimes v_1 &\longmapsto v_1\otimes \tau v_1 \longmapsto 0\\
v_1\otimes v_2 &\longmapsto v_1\otimes \tau v_2 \longmapsto ht\\
v_2\otimes v_1 &\longmapsto v_2\otimes \tau v_1 \longmapsto -ht\\
v_2\otimes v_2 &\longmapsto v_2\otimes \tau v_2 \longmapsto 0,
\end{align*}
and hence factors through the canonical morphism from $E'\otimes E'$ to $\bigwedge^2 E'$.  Denote this induced morphism $\Xi: \bigwedge^2E' \to L$, and define ${\bf F_3}(\xi)=\Xi$.  The inverse morphism, ${\bf G_3}$, is given by pre-composing an element $\Xi\in \Hom(\bigwedge^2E',L)$ with the morphism
\[
\xymatrix{
E'\otimes \tau E' \ar[r]^-{1\otimes \tau} & E'\otimes E' \ar[r]^-{\text{can}} & \bigwedge^2E'.
}
\]
This morphism is defined without reference to a basis, and so is clearly natural.  To check it is an inverse of ${\bf F_3}$, suppose $\Xi\in \Hom(\bigwedge^2E',L)$ is of the form
\[
v_1\wedge v_2 \longmapsto ht.
\]
Then ${\bf G_3}(\Xi)\in \Hom(E'\otimes \tau E', L)$ is given by
\begin{align*}
v_1\otimes \tau v_1 &\longmapsto v_1\otimes v_1 \longmapsto v_1\wedge v_1 = 0 \longmapsto 0\\
v_1\otimes \tau v_2 &\longmapsto v_1\otimes v_2 \longmapsto v_1\wedge v_2 \longmapsto ht\\
v_2\otimes \tau v_1 &\longmapsto v_2\otimes v_1 \longmapsto v_2\wedge v_1 \longmapsto -ht\\
v_2\otimes \tau v_2 &\longmapsto v_2\otimes v_2 \longmapsto v_2\wedge v_2 = 0 \longmapsto 0,
\end{align*}
and hence equals $\xi$.  This proves (iii).

We next prove (i).  Assume $\phi \in M_1$ is in the form of Definition 5.2.  The induced morphism
\[
\xymatrix{
L\otimes E'\otimes E' \ar[r]^-{\phi \otimes 1} & E' \otimes E' \ar[r]^-{\text{can}} & \bigwedge^2 E'
}
\]
maps
\begin{align*}
t\otimes v_1\otimes v_1 &\longmapsto (av_1+bv_2)\otimes v_1 \longmapsto -b\; v_1\wedge v_2\\
t\otimes v_1\otimes v_2 &\longmapsto (av_1+bv_2)\otimes v_2 \longmapsto a\; v_1\wedge v_2\\
t\otimes v_2\otimes v_1 &\longmapsto (cv_1-av_2)\otimes v_1 \longmapsto a\; v_1\wedge v_2\\
t\otimes v_2\otimes v_2 &\longmapsto (cv_1-av_2)\otimes v_2 \longmapsto c\; v_1\wedge v_2,
\end{align*}
and hence factors through the canonical morphism from $L\otimes E'\otimes E'$ to $L\otimes S^2E'$.  Denote this induced morphism $\Phi:L\otimes S^2E' \to \bigwedge^2 E'$, and define ${\bf F_1}(\phi)=\Phi$.

The inverse morphism is constructed as follows.  Observe that we have isomorphisms
\begin{align*}
\Hom(\Hom(L\otimes S^2E', \bigwedge^2E'), &\Hom\left(L\otimes E', E'\right))\\
&\cong \Hom(L^*\otimes (S^2E')^*\otimes \bigwedge^2 E', L^*\otimes E'^*\otimes E')\\
&\cong L\otimes S^2E'\otimes (\bigwedge^2 E')^*\otimes L^*\otimes E'^*\otimes E'\\
&\cong \Hom(L\otimes\bigwedge^2 E'\otimes E', L\otimes E'\otimes S^2E').
\end{align*}
An element of this final group is the morphism ${\bf G_1}$ defined by
\[
l \otimes (e_1\wedge e_2)\otimes e_3\longmapsto l \otimes e_1\otimes e_2e_3 - l \otimes e_2\otimes e_1e_3,
\]
for $l\in L, e_i\in E'$.  In this form, it is clear that ${\bf G_1}$ does not depend on a choice of basis, and is therefore natural.  It remains to check, however, that ${\bf G_1}$ maps $\Hom(L\otimes S^2E', \bigwedge^2E')$ isomorphically onto $M_1$, and that it is the inverse of the map ${\bf F_1}$.

We first trace ${\bf G_1}$ backwards through the chain of isomorphisms. With respect to the basis $\{t, v_1, v_2\}$, the morphism ${\bf G_1}$ is given by
\begin{align*}
t\otimes (v_1\wedge v_2)\otimes v_1&\longmapsto t\otimes v_1\otimes v_1v_2-t\otimes v_2\otimes v_1^2\\
t\otimes (v_1\wedge v_2)\otimes v_2&\longmapsto t\otimes v_1\otimes v_2^2-t\otimes v_2\otimes v_1v_2.
\end{align*}
As an element of $L^*\otimes (\bigwedge^2E')^*\otimes E'^*\otimes L\otimes E' \otimes S^2E'$, this morphism corresponds to
\[
t^*\otimes (v_1\wedge v_2)^*\otimes v_1^*\otimes (t\otimes v_1\otimes v_1v_2-t\otimes v_2\otimes v_1^2)+t^*\otimes (v_1\wedge v_2)^*\otimes v_2^*\otimes (t\otimes v_1\otimes v_2^2-t\otimes v_2\otimes v_1v_2).
\]
As an element of $\Hom(L^*\otimes (S^2E')^*\otimes \bigwedge^2 E', L^*\otimes E'^*\otimes E')$, this element maps
\begin{align*}
t^*\otimes (v_1^2)^* \otimes (v_1\wedge v_2)&\longmapsto -t^*\otimes v_1^*\otimes v_2\\
t^*\otimes (v_1v_2)^*\otimes (v_1\wedge v_2)&\longmapsto t^*\otimes v_1^*\otimes v_1 - t^*\otimes v_2^*\otimes v_2\\
t^*\otimes (v_2^2)^*\otimes (v_1\wedge v_2)&\longmapsto t^*\otimes v_2^*\otimes v_1.
\end{align*}
Now, as an element of $L^*\otimes (S^2E')^*\otimes \bigwedge^2 E'$, the map $\Phi$ corresponds to
\[
(-b (t)^*\otimes (v_1^2)^*+a (t)^*\otimes (v_1v_2)^*+c (t)^*\otimes (v_2^2)^*)\otimes (v_1\wedge v_2),
\]
and hence the image of $\Phi$ under ${\bf G_1}$ is
\[
bt^*\otimes v_1^* \otimes v_2 + a(t^*\otimes v_1^*\otimes v_1 - t^*\otimes v_2^*\otimes v_2)+c t^*\otimes v_2^*\otimes v_1
\]
in $L^*\otimes E'^*\otimes E'$.  As an element of $\Hom(L\otimes E', E')$, this element maps
\begin{align*}
t\otimes v_1&\longmapsto av_1+bv_2\\
t\otimes v_2&\longmapsto cv_1-av_2,
\end{align*}
and hence equals $\phi$.
\end{proof}

\tpoint{Definition} Let $S_3\text{Cov}_B(A)\subset \Hom_B(S^2A,A)$ denote the submodule of morphisms defining commutative, associative $B$-algebra structures on $A$ compatible with the given $S_3$-action (and hence inducing $S_3$-covers $\pi: \Spec{A}\to \Spec{B}$). Let $S_3\text{Cov}_B(A)^0$ denote those which define integral such algebras.

Also define
\[
\text{Build}_B(A)=\Hom(L\otimes S^2E', \bigwedge^2E') \oplus \Hom(S^3E', \bigwedge^2E')\oplus \Hom(\bigwedge^2 E', L).
\]

By the previous lemma and Corollary \ref{localtoglobal}, when then have the following:

\tpoint{Corollary}{\em   There exists a natural morphism ${\bf F}:S_3\text{Cov}_B(A)^0\to \text{Build}_B(A)$.
}

As we'll soon see, the morphism ${\bf F}$ extracts from an $S_3$-cover the ``building data" necessary to reconstruct the cover.

Note that any triple $(\Phi, \Psi, \Xi)\in \text{Build}_B(A)$ is of the form
\begin{align*}
\Phi(t\otimes v_1^2)&= A (v_1\wedge v_2)\\
\Phi(t\otimes v_1v_2)&= B (v_1\wedge v_2)\\
\Phi(t\otimes v_2^2)&= C (v_1\wedge v_2)\\
\Psi(v_1^3)&= D (v_1\wedge v_2)\\
\Psi(v_1^2v_2)&= E (v_1\wedge v_2)\\
\Psi(v_1v_2^2)&= F (v_1\wedge v_2)\\
\Psi(v_2^3)&= G (v_1\wedge v_2)\\
\Xi(v_1\wedge v_2)&= ht
\end{align*}
with respect to the basis $\{t, v_1, v_2\}$ for $L\oplus E'$.  We use capital letters here to avoid confusion with our notation for a triple ($\phi, \psi, \xi)\in \Hom(L\otimes E',E')\oplus \Hom(S^2E',E')\oplus \Hom(E'\otimes \tau E',L)$. (From the context, it should be clear when $B$ refers to a constant and when it refers to the ring $\Oh_X(U)$.) Under the natural isomorphisms of Lemma \ref{basis-free-description}, if we let $(\phi,\psi,\xi)={\bf F}^{-1}(\Phi,\Psi,\Xi)$, we have
\begin{align*}
\phi(t\otimes v_1)&= Bv_1-Av_2\\
\phi(t\otimes v_2)&= Cv_1-Bv_2\\
\psi(v_1^2)&= Ev_1-Dv_2\\
\psi(v_1v_2)&= Fv_1-Ev_2\\
\psi(v_2^2)&= Gv_1-Fv_2\\
\xi(v_1\otimes \tau v_2)&= ht.
\end{align*}
The correspondence with our earlier notation is thus $A=-b, B=a, C=c, D=-e, E=d, F=-g,G=f$.  Using this dictionary, the three conditions of Corollary \ref{localtoglobal} become the following:
\begin{itemize}
\item[(i)] $AF-2BE+CD=0$;
\item[(ii)] $AG-2BF+CE=0$;
\item[(iii)] $h(B^2-AC)=\frac{3}{2}(B(EF-DG)-A(F^2-EG)+C(DF-E^2)).$
\end{itemize}

\tpoint{Lemma}{\em There exists a natural morphism, ${\bf A_1}$, from $\Hom\left(L\otimes S^2E', \bigwedge^2E'\right)\oplus \Hom\left(S^3E',\bigwedge^2E'\right)$ to $\Hom\left(L\otimes \left(\bigwedge^2E'\right)^{\otimes 2}\otimes E',\left(\bigwedge^2E'\right)^{\otimes 2}\right)$
whose kernel consists of precisely those pairs $(\Phi, \Psi)$ satisfying conditions (i) and (ii).}

\begin{proof}
First, consider the morphism $f_1:(\bigwedge^2 E')^{\otimes 2}\to (S^2E')^{\otimes 2}$ defined by
\[
(e_1\wedge e_2)\otimes (e_3\wedge e_4)\longmapsto e_1e_3\otimes e_2e_4-e_1e_4\otimes e_2e_3-e_2e_3\otimes e_1e_4+e_2e_4\otimes e_1e_3,
\]
and the canonical morphism $f_2:S^2E'\otimes E'\to S^3E'$.  Given any pair $(\Phi, \Psi)\in \Hom(L\otimes S^2E', \bigwedge^2E')\oplus \Hom(S^3E',\bigwedge^2E')$, we then have an induced morphism
\[
\xymatrix{
L\otimes (\bigwedge^2 E')^{\otimes 2}\otimes E' \ar[r]^-{1\otimes f_1\otimes 1} & L \otimes (S^2E')^{\otimes 2}\otimes E' \ar[r]^-{1\otimes 1\otimes f_2} & L\otimes S^2E' \otimes S^3E' \ar[r]^-{\Phi \otimes \Psi} & (\bigwedge^2 E')^{\otimes 2}
}
\]
Let ${\bf A_1}(\Phi,\Psi)$ denote this morphism. With respect to the basis, this morphism is defined by
\begin{align*}
t\otimes (v_1\wedge v_2)\otimes (v_1\wedge v_2)\otimes v_1 &\longmapsto t\otimes (v_1^2\otimes v_2^2-2(v_1v_2\otimes v_1v_2)+v_2^2\otimes v_1^2)\otimes v_1\\
& \longmapsto t\otimes (v_1^2\otimes v_1v_2^2-2(v_1v_2\otimes v_1^2v_2)+v_2^2\otimes v_1^3)\\
&\longmapsto (AF-2BE+CD)(v_1\wedge v_2)^{\otimes 2}\\
t\otimes (v_1\wedge v_2)\otimes (v_1\wedge v_2)\otimes v_2 &\longmapsto t\otimes (v_1^2\otimes v_2^2-2(v_1v_2\otimes v_1v_2)+v_2^2\otimes v_1^2)\otimes v_2\\
& \longmapsto t\otimes (v_1^2\otimes v_2^3-2(v_1v_2\otimes v_1v_2^2)+v_2^2\otimes v_1^2v_2)\\
&\longmapsto (AG-2BF+CE)(v_1\wedge v_2)^{\otimes 2}
\end{align*}
Thus, conditions (i) and (ii) together are equivalent to the morphism ${\bf A_1}(\Phi, \Psi)$ being the zero map.
\end{proof}

For the proof of the following lemma, note that given any pair $(\Phi, \Psi)\in \Hom(L\otimes S^2E', \bigwedge^2E')\oplus \Hom(S^3E',\bigwedge^2E')$, we naturally have induced maps $\bigwedge^2(\phi),\bigwedge^2(\psi)$.  In terms of the basis, these are given by
\begin{align*}
(t\otimes v_1)\wedge (t\otimes v_2)&\longmapsto -(B^2-AC)(v_1\wedge v_2)\\
\intertext{and}
v_1^2\wedge v_1v_2 &\longmapsto (DF-E^2)(v_1\wedge v_2)\\
v_1^2\wedge v_2^2 &\longmapsto (DG-EF)(v_1\wedge v_2)\\
v_1v_2\wedge v_2^2 &\longmapsto (EG-F^2)(v_1\wedge v_2),
\end{align*}
respectively.

\tpoint{Lemma}{\em There exists a natural morphism, ${\bf A_2}$, from $\text{Build}_B (A)$ to $\Hom (L\otimes\left(\bigwedge^2E'\right)^{\otimes 3},$ $\left(\bigwedge^2 E'\right)^{\otimes 2})$
whose kernel consists of precisely those triples $(\Phi, \Psi, \Xi)$ satisfying condition (iii).}

\begin{proof}
Consider the two natural morphisms, $f_1:(\bigwedge^2E')^{\otimes 3}\to \bigwedge^3(S^2E')$ and $f_2:\bigwedge^3(S^2E')\to S^2E' \otimes \bigwedge^2(S^2E')$, defined by
\[
(e_1\wedge e_2)\otimes (e_3\wedge e_4)\otimes (e_5\wedge e_6)\longmapsto e_1^2\wedge e_2^2\wedge (e_3e_5-e_3e_6-e_4e_5+e_4e_6),
\]
and
\[
e_1e_2\wedge e_3e_4\wedge e_5e_6 \longmapsto \frac{3}{4}\left(e_1e_2\otimes (e_3e_4\wedge e_5e_6)-e_3e_4\otimes (e_1e_2\wedge e_5e_6)+e_5e_6\otimes (e_1e_2\wedge e_3e_4)\right),
\]
respectively.  Given any pair $(\Phi, \Psi)\in \Hom(L\otimes S^2E', \bigwedge^2E')\oplus \Hom(S^3E',\bigwedge^2E')$ we then have an induced morphism
\[
\xymatrix{
L\otimes (\bigwedge^2 E')^{\otimes 3} \ar[r]^-{1\otimes f_1} & L\otimes \bigwedge^3(S^2E') \ar[r]^-{1\otimes f_2} & L\otimes S^2E' \otimes \bigwedge^2(S^2E') \ar[r]^-{\Phi\otimes \bigwedge^2(\psi)} & (\bigwedge^2 E')^{\otimes 2}
}.
\]
Denote this morphism by ${\bf A_{2,1}}(\Phi, \Psi)$.  In terms of the basis, this morphism is given by
\begin{align*}
t\otimes (v_1\wedge v_2)^{\otimes 3} &\longmapsto t\otimes (v_1^2\wedge v_2^2 \wedge (v_1^2-v_1v_2-v_2v_1+v_2^2) = 2 t\otimes (v_1^2\wedge v_1v_2 \wedge v_2^2)\\
&\longmapsto \frac{3}{2} t\otimes \left(v_1^2\otimes (v_1v_2\wedge v_2^2)-v_1v_2\otimes (v_1^2\wedge v_2^2)+v_2^2\otimes (v_1^2\wedge v_1v_2)\right)\\
&\longmapsto \frac{3}{2}\left(A(EG-F^2)-B(DG-EF)+C(DF-E^2)\right)(v_1\wedge v_2)^{\otimes 2}\\
&\quad = \frac{3}{2}\left(B(EF-DG)-A(F^2-EG)+C(DF-E^2)\right)(v_1\wedge v_2)^{\otimes 2}.
\end{align*}

Given any pair $(\Phi, \Xi)\in \Hom(L\otimes S^2E', \bigwedge^2E')\oplus\Hom(\bigwedge^2 E', L)$, we have an induced morphism
\[
\xymatrix{
\bigwedge^2E' \otimes E' \otimes L\otimes E' \ar[r]^-{\Xi\otimes 1^{\otimes 3}} & L\otimes E'\otimes L\otimes E' \ar[r]^-{\text{can}} & \bigwedge^2(L\otimes E') \ar[r]^-{\bigwedge^2(\phi)} & \bigwedge^2 E'.
}
\]
In terms of the basis, this morphism maps
\begin{align*}
t\otimes(v_1\wedge v_2)\otimes v_1\otimes v_1&\longmapsto ht\otimes v_1\otimes t\otimes v_1\longmapsto h(t\otimes v_1)\wedge (t\otimes v_1) = 0\longmapsto 0\\
t\otimes(v_1\wedge v_2)\otimes v_1\otimes v_2&\longmapsto ht\otimes v_1\otimes t\otimes v_2\longmapsto h(t\otimes v_1)\wedge (t\otimes v_2)\longmapsto -h(B^2-AC)(v_1\wedge v_2)\\
t\otimes(v_1\wedge v_2)\otimes v_2\otimes v_1&\longmapsto ht\otimes v_2\otimes t\otimes v_1\longmapsto h(t\otimes v_2)\wedge (t\otimes v_1)\longmapsto h(B^2-AC)(v_1\wedge v_2)\\
t\otimes(v_1\wedge v_2)\otimes v_2\otimes v_2&\longmapsto ht\otimes v_2\otimes t\otimes v_2\longmapsto h(t\otimes v_2)\wedge (t\otimes v_2) = 0\longmapsto 0,
\end{align*}
and hence factors through the canonical morphism (induced from) $E'\otimes E' \to \bigwedge^2E'$.  This gives a morphism
\[
{\bf A'_{2,2}}(\Phi,\Xi):L\otimes (\bigwedge^2E')^{\otimes 2} \to \bigwedge^2E'.
\]
Let ${\bf A_{2,2}}(\Phi, \Xi)$ then denote the induced morphism
\[
\xymatrix{
L\otimes (\bigwedge^2E')^{\otimes 3} \ar[r]^-{{\bf A'_{2,2}}\otimes 1} & (\bigwedge^2E')^{\otimes 2}.
}
\]
In terms of the basis, this morphism is given by
\[
t\otimes (v_1\wedge v_2)^{\otimes 3}\longmapsto -h(B^2-AC)(v_1\wedge v_2)^{\otimes 2}.
\]

Lastly, define ${\bf A_2}(\Phi, \Psi, \Xi)\in \Hom(L\otimes(\bigwedge^2E')^{\otimes 3}, (\bigwedge^2 E')^{\otimes 2})$ as the sum ${\bf A_{2,1}}(\Phi, \Psi)+{\bf A_{2,2}}(\Phi, \Xi)$.  In terms of the basis, this morphism is given by
\[
t\otimes (v_1\wedge v_2)^{\otimes 3}\longmapsto \left(\frac{3}{2}\left(B(EF-DG)-A(F^2-EG)+C(DF-E^2)\right)-h(B^2-AC)\right)(v_1\wedge v_2)^{\otimes 2},
\]
and hence condition (iii) is equivalent to the vanishing of ${\bf A_2}(\Phi, \Psi, \Xi)$.
\end{proof}

\tpoint{Definition} For notational simplicity, let us define
\[
\text{Compat}_B(A)=\Hom(L\otimes (\bigwedge^2 E')^{\otimes 2}\otimes E',(\bigwedge^2 E')^{\otimes 2})\oplus \Hom(L\otimes (\bigwedge^2 E')^{\otimes 3}, (\bigwedge^2 E')^{\otimes 2}).
\]

The previous two lemmas then give the following:

\tpoint{Corollary}{\em There is a natural morphism ${\bf A}:\text{Build}_B(A)\to \text{Compat}_B(A)$ whose kernel consists of precisely the morphisms satisfying conditions (i)-(iii).}

In other words, the morphism ${\bf A}$ tests the building data for the compatibility conditions (arising from associativity constraints) necessary for the data to induce an $S_3$-cover.

We next recover from any building data the three morphisms defining an algebra structure:
\[
\alpha: S^2L\to B, \quad \beta:L\otimes E\to E, \quad \gamma: S^2E\to A.
\]

\tpoint{Lemma}{\em There exists a natural morphism, ${\bf B_1}$, from $\Hom\left(L\otimes S^2E', \bigwedge^2E'\right)$ to  $\Hom\left(S^2L,B\right)$,
taking an element $\Phi$ of the form
\begin{align*}
\Phi(t\otimes v_1^2)&= A (v_1\wedge v_2)\\
\Phi(t\otimes v_1v_2)&= B (v_1\wedge v_2)\\
\Phi(t\otimes v_2^2)&= C (v_1\wedge v_2)
\end{align*}
to an element $\alpha$ of the form
\[
\alpha(t^2)=-3(B^2-AC).
\]}

\begin{proof}
Observe that we have isomorphisms
\begin{align*}
\Hom\left(\Hom\left(\bigwedge^2(L\otimes E'), \bigwedge^2E'\right), \Hom(S^2L,B)\right)&\cong \Hom\left(\left(\bigwedge^2(L\otimes E')\right)^*\otimes \bigwedge^2 E', (S^2L)^*\right)\\
&\cong \Hom\left(\bigwedge^2E'\otimes S^2L, \bigwedge^2(L\otimes E')\right).
\end{align*}
An element of this last group is the morphism ${\bf \tilde{B}_1}$ defined by
\[
(e_1\wedge e_2)\otimes l_1l_2 \longmapsto \frac{3}{2}\left((l_1\otimes e_1)\wedge (l_2\otimes e_2)-(l_1\otimes e_2)\wedge (l_2\otimes e_1)\right).
\]
In terms of the basis, this morphism is defined by
\[
(v_1\wedge v_2)\otimes t^2 \longmapsto \frac{3}{2}\left((t\otimes v_1)\wedge (t\otimes v_2)-(t\otimes v_2)\wedge (t\otimes v_1)\right)=3 (t\otimes v_1)\wedge (t\otimes v_2).
\]
This morphism, considered as an element of $\Hom\left((\bigwedge^2(L\otimes E'))^*\otimes \bigwedge^2 E', (S^2L)^*\right)$, is given by
\[
\left( (t\otimes v_1)\wedge (t\otimes v_2)\right)^*\otimes (v_1\wedge v_2)\longmapsto 3(t^2)^*.
\]
The morphism $\bigwedge^2 (\phi)$, when considered as an element of $(\bigwedge^2(L\otimes E'))^*\otimes \bigwedge^2 E'$, is
\[
-(B^2-AC)\left( (t\otimes v_1)\wedge (t\otimes v_2)\right)^*\otimes (v_1\wedge v_2).
\]
Thus, ${\bf \tilde{B}_1}$ maps $\bigwedge^2 (\phi)$ to the element
\[
-3(B^2-AC)(t^2)^*,
\]
which corresponds to the map
\[
t^2\longmapsto -3(B^2-AC).
\]
Thus, the composition ${\bf B_1}={\bf \tilde{B}_1}\circ \bigwedge^2\circ {\bf F_1}^{-1}$ is the desired natural transformation.
\end{proof}

\tpoint{Lemma}{\em There exists a natural morphism, ${\bf B_2}$, from $\Hom\left(L\otimes S^2E', \bigwedge^2E'\right)$ to  $\Hom(L\otimes E, E)$,
taking an element $\Phi$ of the form
\begin{align*}
\Phi(t\otimes v_1^2)&= A (v_1\wedge v_2)\\
\Phi(t\otimes v_1v_2)&= B (v_1\wedge v_2)\\
\Phi(t\otimes v_2^2)&= C (v_1\wedge v_2)
\end{align*}
to an element $\beta$ of the form
\begin{align*}
\beta(t\otimes v_1)&= Bv_1-Av_2-2B \tau v_1 +2A\tau v_2\\
\beta(t\otimes v_2)&= Cv_1-Bv_2-2C\tau v_1 +2B\tau v_2\\
\beta(t\otimes \tau v_1)&= 2Bv_1-2Av_2-B\tau v_1 +A\tau v_2\\
\beta(t\otimes \tau v_2)&= 2Cv_1-2Bv_2-C\tau v_1 +B\tau v_2.
\end{align*}}

\begin{proof}
Observe that $\phi={\bf F_1}^{-1}(\Phi)$ induces a morphism
\[
\xymatrix{
L\otimes E' \ar[r]^-{\phi} & E' \ar[r]^-{\tau} & \tau E' \subset E,
}
\]
and hence two morphisms:
\[
\xymatrix{
\beta_1: L\otimes E' \ar[rr]^-{(\phi, -2(\tau \circ \phi))} & & E'\oplus \tau E' = E,
}
\]
and
\[
\xymatrix{
\beta_2: L\otimes \tau E' \ar[r]^-{\tau \otimes \tau} & L\otimes E' \ar[rr]^-{(\phi, -2(\tau \circ \phi))} & & E'\oplus \tau E' = E \ar[r]^-{\tau} & E.
}
\]
Let ${\bf B_2}(\Phi)= \langle \beta_1, \beta_2\rangle : L\otimes E \cong (L\otimes E')\oplus (L\otimes \tau E')\to E$. In terms of the basis, this morphism is given by
\begin{align*}
t\otimes v_1&\longmapsto  Bv_1-Av_2-2B \tau v_1 +2A\tau v_2\\
t\otimes v_2&\longmapsto Cv_1-Bv_2-2C\tau v_1 +2B\tau v_2\\
t\otimes \tau v_1&\longmapsto 2Bv_1-2Av_2-B\tau v_1 +A\tau v_2\\
t\otimes \tau v_2&\longmapsto 2Cv_1-2Bv_2-C\tau v_1 +B\tau v_2.
\end{align*}
\end{proof}

\tpoint{Lemma}{\em There exists a natural morphism, ${\bf B_3}$, from $\Hom\left(S^3E', \bigwedge^2E'\right)\oplus \Hom\left(\bigwedge^2 E', L\right)$ to $\Hom(S^2E,A)$,
taking a pair $(\Psi, \Xi)$ of the form
\begin{align*}
\Psi(v_1^3)&= D (v_1\wedge v_2)\\
\Psi(v_1^2v_2)&= E (v_1\wedge v_2)\\
\Psi(v_1v_2^2)&= F (v_1\wedge v_2)\\
\Psi(v_2^3)&= G (v_1\wedge v_2)\\
\Xi(v_1\wedge v_2)&= ht
\end{align*}
to an element $\gamma$ of the form
\begin{align*}
\gamma(v_1^2)&= 6(E^2-DF)+Ev_1-Dv_2-2E\tau v_1 +2D\tau v_2\\
\gamma(v_1v_2)&=3(EF-DG)+Fv_1-Ev_2-2F\tau v_1 + 2E\tau v_2\\
\gamma(v_2^2)&= 6(F^2-EG)+Gv_1-Fv_2-2G\tau v_1 +2F\tau v_2\\
\gamma(v_1\cdot \tau v_1)&= 3(E^2-DF)-Ev_1+Dv_2-E\tau v_1+D\tau v_2\\
\gamma(v_1\cdot \tau v_2)&= \frac{3}{2}(EF-DG)+ht-Fv_1+Ev_2-F\tau v_1 +E\tau v_2\\
\gamma(v_2\cdot \tau v_1)&= \frac{3}{2}(EF-DG)-ht-Fv_1+Ev_2-F\tau v_1 +E\tau v_2\\
\gamma(v_2\cdot \tau v_2)&= 3(F^2-EG)-Gv_1+Fv_2-G\tau v_1 +F\tau v_2\\
\gamma((\tau v_1)^2)&= 6(E^2-DF)-2Ev_1+2Dv_2+E\tau v_1-D\tau v_2\\
\gamma((\tau v_1)(\tau v_2))&= 3(EF-DG)-2Fv_1+2Ev_2+F\tau v_1-E\tau v_2\\
\gamma((\tau v_2)^2)&= 6(F^2-EG)-2Gv_1+2Fv_2+G\tau v_1 -F\tau v_2.
\end{align*}}

\begin{proof}
Observe that we have isomorphisms
\begin{align*}
\Hom\left(\Hom\left(\bigwedge^2(S^2E'), \bigwedge^2 E'\right),\Hom(S^2E', B)\right)&\cong \Hom\left(\left(\bigwedge^2(S^2E')\right)^*\otimes \bigwedge^2E', (S^2E')^*\right)\\
&\cong \Hom\left(\left(\bigwedge^2E'\right)\otimes S^2E', \bigwedge^2(S^2E')\right).
\end{align*}
A natural element of this last group is the morphism ${\bf \tilde{B}_{3,1}}$ defined by
\[
(e_1\wedge e_2)\otimes e_3e_4 \longmapsto -3(e_1e_3\wedge e_2e_4+e_1e_4\wedge e_2e_3).
\]
In terms of the basis, this morphism is given by
\begin{align*}
(v_1\wedge v_2)\otimes v_1^2 &\longmapsto -3(v_1^2\wedge v_1v_2+v_1^2\wedge v_1v_2)=-6 (v_1^2\wedge v_1v_2)\\
(v_1\wedge v_2)\otimes v_1v_2 &\longmapsto -3(v_1^2\wedge v_2^2+v_1v_2\wedge v_1v_2)= -3(v_1^2\wedge v_2^2)\\
(v_1\wedge v_2)\otimes v_2^2 &\longmapsto -3(v_1v_2\wedge v_2^2 +v_1v_2\wedge v_2^2)= -6 (v_1v_2\wedge v_2^2).
\end{align*}
As an element of $\Hom((\bigwedge^2(S^2E'))^*\otimes \bigwedge^2E', (S^2E')^*)$, this corresponds to the map
\begin{align*}
(v_1^2\wedge v_1v_2)^*\otimes (v_1\wedge v_2)&\longmapsto -6(v_1^2)^*\\
(v_1^2\wedge v_2^2)^*\otimes (v_1\wedge v_2)&\longmapsto -3(v_1v_2)^*\\
(v_1v_2\wedge v_2^2)^*\otimes (v_1\wedge v_2)&\longmapsto -6(v_2^2)^*.
\end{align*}
The morphism $\bigwedge^2(\psi)$, as an element of $(\bigwedge^2(S^2E'))^*\otimes \bigwedge^2E'$,  is given by
\[
(DF-E^2)(v_1^2\wedge v_1v_2)^*\otimes (v_1\wedge v_2)+(DG-EF)(v_1^2\wedge v_2^2)^*\otimes (v_1\wedge v_2)+(EG-F^2)(v_1v_2\wedge v_2^2)^*\otimes (v_1\wedge v_2).
\]
Thus, ${\bf \tilde{B}_{3,1}}$ maps $\bigwedge^2(\psi)$ to the element
\[
6(E^2-DF)(v_1^2)^*+3(EF-DG)(v_1v_2)^*+6(F^2-EG)(v2^2)^*,
\]
which corresponds to the morphism $\tilde{\gamma}_1$ defined by
\begin{align*}
v_1^2&\longmapsto 6(E^2-DF)\\
v_1v_2&\longmapsto 3(EF-DG)\\
v_2^2&\longmapsto 6(F^2-EG).
\end{align*}
This induces morphisms
\[
\xymatrix{
\gamma_{1,1}': E'\otimes E' \ar[rr]^-{\text{can}} & & S^2E' \ar[r]^-{\tilde{\gamma}_1} & B\\
\gamma_{1,2}': E'\otimes \tau E' \ar[r]^-{1\otimes \tau} & E' \otimes E' \ar[r]^-{\text{can}} & S^2E' \ar[r]^-{\frac{1}{2}\tilde{\gamma}_1} & B\\
\gamma_{1,3}': \tau E'\otimes E' \ar[r]^-{\tau \otimes 1} & E' \otimes E' \ar[r]^-{\text{can}} & S^2E' \ar[r]^-{\frac{1}{2}\tilde{\gamma}_1} & B\\
\gamma_{1,4}': \tau E'\otimes \tau E' \ar[r]^-{\tau \otimes \tau} & E' \otimes E' \ar[r]^-{\text{can}} & S^2E' \ar[r]^-{\tilde{\gamma}_1} & B,
}
\]
which together define a morphism $\gamma_1':E\otimes E \cong (E'\otimes E')\oplus (E'\otimes \tau E') \oplus (\tau E'\otimes E') \oplus (\tau E'\otimes \tau E') \to B$.  By construction, this morphism factors through the canonical morphism from $E\otimes E$ to $S^2E$, and gives the first-coordinate morphism $\gamma_1:S^2E\to B$.  Let ${\bf B_{3,1}}(\Psi)=\gamma_1$.

Next, observe that $\xi={\bf F_3}^{-1}(\Xi)$ induces a morphism
\[
\xymatrix{
\xi':\tau E'\otimes E' \ar[r]^-{\text{can}} & E'\otimes \tau E' \ar[r]^-{\xi} & L.
}
\]
Define $\gamma_2':E\otimes E\to L$ by $\gamma_2'=\langle 0, \xi, \xi', 0\rangle$.  By construction, this morphism factors through the canonical morphism from $E\otimes E$ to $S^2E$, and gives the second-coordinate morphism $\gamma_2:S^2E\to L$. Let ${\bf B_{3,2}}(\Xi)=\gamma_2$.

Lastly, observe that the morphism $\psi={\bf F_2}^{-1}(\Psi)$ induces morphisms
\[
\xymatrix{
\gamma_{3,1}': E'\otimes E' \ar[rr]^-{\text{can}} & & S^2E' \ar[r]^-{\psi-2(\tau \circ \psi)} & E'\oplus \tau E'\\
\gamma_{3,2}': E'\otimes \tau E' \ar[r]^-{1\otimes \tau} & E' \otimes E' \ar[r]^-{\text{can}} & S^2E' \ar[r]^-{-\psi-(\tau \circ \psi)} & E'\oplus \tau E'\\
\gamma_{3,3}': \tau E'\otimes E' \ar[r]^-{\tau \otimes 1} & E' \otimes E' \ar[r]^-{\text{can}} & S^2E' \ar[r]^-{-\psi-(\tau \circ \psi)} & E'\oplus \tau E'\\
\gamma_{3,4}': \tau E'\otimes \tau E' \ar[r]^-{\tau \otimes \tau} & E' \otimes E' \ar[r]^-{\text{can}} & S^2E' \ar[r]^-{-2\psi+(\tau \circ \psi)} & E'\oplus \tau E'.
}
\]
These together define a morphism $\gamma_3':E\otimes E\to E$ by $\gamma_3'=\langle \gamma_{3,1}', \gamma_{3,2}', \gamma_{3,3}', \gamma_{3,4}'\rangle$.  By construction, this factors through the canonical morphism from $E\otimes E$ to $S^2E$, and gives the third-coordinate morphism $\gamma_3:S^2E\to L$. Let ${\bf B_{3,3}}(\Psi)=\gamma_3$.
\end{proof}

\tpoint{Corollary}{\em There is a natural morphism ${\bf B}: \text{Build}_B(A)\to \Hom_B(S^2A,A)$ which maps $\ker{{\bf A}}\subset \text{Build}_B(A)$ into $S_3\text{Cov}_B(A)$.}

In other words, given compatible building data, ${\bf B}$ builds the associated $S_3$-cover.

By naturality, the results above sheafify to give the following:

\tpoint{Main Theorem}\label{maintheorem}{\em Let $X$ be an integral, Noetherian $R$-scheme.  Suppose $\mathcal{L}$ is an invertible $\Oh_X$-module on which $S_3$ acts via the sign character, and $\E$ is a locally free $\Oh_X$-module of rank $4$ together with an $S_3$-action such that (under the induced $R[S_3]$-action) $e_i$ acts as $\delta_{3i}\text{Id}_{\E}$.  Let $\E=\E'\oplus \tau \E'$ be the $\Oh_X$-module decomposition of $\E$ induced by $e_3=e_{31}+e_{32}$, and let $\A=\Oh_X\oplus \mathcal{L}\oplus \E$.

There exist morphisms ${\bf F}:S_3\text{Cov}_X(\A)^0\to \text{Build}_X(\A)$, ${\bf B}:\text{Build}_X(\A) \to \Hom_X(S^2\A,\A)$, and ${\bf A}:\text{Build}_X(\A)\to \text{Compat}_X(\A)$ with the following properties:
\begin{itemize}
\item[(i)] The morphism ${\bf F}$ factors through $\ker{{\bf A}}$; i.e., if $\A$ is endowed with an algebra structure (compatible with the given $S_3$-action) such that $\pi:\RelSpec{X}{\A}\to X$ is a flat $S_3$-cover of integral, Noetherian $R$-schemes, then the algebra structure on $\A$ defines an element $(\alpha, \beta, \gamma)\in S_3\text{Cov}_X(\A)^0$ with ${\bf F}(\alpha, \beta,\gamma)\in \ker{{\bf A}}$; and
\item[(ii)] The morphism ${\bf B}|_{\ker{{\bf A}}}$ factors through $S_3\text{Cov}_X(\A)$; i.e., given any element $(\Phi, \Psi, \Xi)\in \ker{{\bf A}}$, the element ${\bf B}(\Phi, \Psi, \Xi)$ defines a commutative, associative $\Oh_X$-algebra structure on $\A$ compatible with the given $S_3$-action, and hence defines a flat $S_3$-cover $\pi:\RelSpec{X}{\A}\to X$.
\end{itemize}
}

\bpoint{Closing remarks} This result represents only the starting point of the study of $S_3$-covers.  In analogy with the abelian case, an immediate question is under what circumstances the morphisms ${\bf F}$ and ${\bf B}$ are mutually inverse; that is, when does the building data uniquely determine the cover?  Other points of study include (but are certainly not limited to) the global description of the branch locus, the location and nature of singularities of the covering scheme, and the relationship between the (cohomological) invariants of the base scheme and those of the covering scheme.

} 

\bibliographystyle{plain}
\bibliography{S_3-covers}

\begin{thebibliography}{10}

\bibitem{bogomolov}
F.A. Bogomolov.
\newblock Holomorphic tensors and vector bundles on projective varieties.
\newblock {\em Math. USSR, Izv.}, (13):499--555, 1979.

\bibitem{thesis}
R.~Easton.
\newblock ${S}_3$-covers of schemes.
\newblock Thesis, 2007.

\bibitem{easton}
R.~Easton.
\newblock Surfaces violating {B}ogomolov-{M}iyaoka-{Y}au in positive
  characteristic.
\newblock {\em Proc. Amer. Math. Soc.}, 136:2271--2278, 2008.

\bibitem{eastonvakil}
R.~Easton and R.~Vakil.
\newblock Absolute {G}alois acts faithfully on the components of the moduli
  space of surfaces: A {B}elyi-type theorem in higher dimension.
\newblock {\em Int. Math. Res. Notices}, 2007, 2007.
\newblock Article ID rnm080, 10 pages, doi:10.1093/imrn/rnm080.

\bibitem{hirzebruch}
F.~Hirzebruch.
\newblock Arrangements of lines and algebraic surfaces.
\newblock In {\em Arithmetic and Geometry, Vol. II}, number~36, pages 113--140.
  Progr. Math., 1983.

\bibitem{lang}
W.E. Lang.
\newblock Examples of surfaces of general type with vector fields.
\newblock In {\em Arithmetic and Geometry, Vol. II}, number~36, pages 167--173.
  Progr. Math., 1983.

\bibitem{miranda}
R.~Miranda.
\newblock Triple covers in algebraic geometry.
\newblock {\em Amer. Jour. of Math.}, 107(5):1123--1158, 1985.

\bibitem{miyaoka}
Y.~Miyaoka.
\newblock On the {C}hern numbers of surfaces of general type.
\newblock {\em Invent. Math.}, (42):225--237, 1977.

\bibitem{pardini}
R.~Pardini.
\newblock Abelian covers of algebraic varieties.
\newblock {\em J. Reine Angew. Math.}, (417):191 -- 213, 1991.

\bibitem{yau}
S.-T. Yau.
\newblock The {C}alabi conjecture and some new results in algebraic geometry.
\newblock {\em Proc. Natl. Acad. Sci. USA}, (74):1798--1799, 1977.

\end{thebibliography}



\end{document}